\documentclass[leqno]{elsart}
\usepackage{amsmath,amscd, amsfonts}

\newtheorem{Theorem}{Theorem}[section]
\newtheorem{Lemma}[Theorem]{Lemma}
\newtheorem{Corollary}[Theorem]{Corollary}
\newtheorem{Proposition}[Theorem]{Proposition}

\def\diam{\operatorname{diam}}
\def\Max{\operatorname{Max}}
\def\deg{\operatorname{deg}}

\def\Hom{\operatorname{Hom}}
\def\bg{\operatorname{beg}}
\def\ed{\operatorname{end}}

\def\Nset{{\mathbb N}}
\def\Zset{{\mathbb Z}}

\def\Pset{{\mathbb P}}
\def\Acal{{\mathcal A}}
\def\Lcal{{\mathcal L}}
\def\Rcal{{\mathcal R}}

\def\mfr{{\mathfrak m}}
\def\alphaun{\underline{\alpha}}
\def\betaun{\underline{\beta}}
\def\xund{\underline{x}}
\def\yund{\underline{y}}

\begin{document}
\begin{frontmatter}

\title{ On local cohomology  of a tetrahedral curve}
\thanks{Both authors were supported by NAFOSTED
(Vietnam)}
\author{D\^o Ho\`ang Giang} and
\ead{dhgiang@math.ac.vn}
\author{L\^e Tu\^an Hoa}
\address{Institute of Mathematics Hanoi, 18 Hoang Quoc Viet Road, 10307 Hanoi,
Vietnam}
\ead{lthoa@math.ac.vn}
\begin{abstract}  It is shown that the diameter $\diam (H^1_\mfr(R/I))$ of the first
local cohomology module of a tetrahedral curve $C= C(a_1,...,a_6)$
can be explicitly expressed in terms of the $a_i$ and is the
smallest non-negative integer $k$ such that $\mfr^k
H^1_\mfr(R/I)=0$. From that one can describe all arithmetically
Cohen-Macaulay or Buchsbaum tetrahedral curves.
\end{abstract}

 \maketitle

\end{frontmatter}

\section*{Introduction}

A tetrahedral curve $C= C(a_1,...,a_6)$ is a curve in $\Pset^3$ defined by the  ideal
$$I =(x_1,x_2)^{a_1}\cap (x_1,x_3)^{a_2}\cap (x_1,x_4)^{a_3}\cap  
(x_2,x_3)^{a_4}\cap (x_2,x_4)^{a_5}\cap (x_3,x_4)^{a_6} $$
of the polynomial ring $R=K[x_1,x_2,x_3,x_4]$ over a field $K$,
where $a_1,...,a_6$ are non-negative integers and not all of  them
are zero. The case $a_2=a_5=0$ was first considered by Schwartau
\cite{S}. He gave  a characterization of the Cohen-Macaulay
property of $C$ in terms of  $a_1,a_3,a_4,a_6$. The general case
of  tetrahedral curves, when $a_2$ and $a_5$ are not necessarily
zero,  was  introduced in \cite {MN}. Using basic double linkage,
Migliore and Nagel gave there an efficient numerical algorithm for
determining when a particular tetrahedral curve is arithmetically
Cohen-Macaulay and asked for an explicit characterization in terms
of $a_1,...,a_6$.  This problem was solved later by Francisco  in
\cite{F}.  Moreover,  it was  shown in the papers \cite{MN,FMN} that
these curves have  many nice properties.

In this paper we study the structure of the first local cohomology module
$H^1_\mfr(R/I)$ with the support in the maximal  homogeneous ideal
$\mfr = (x_1,x_2,x_3,x_4)$.  This study is important because we can characterize many properties, 
such as the Cohen-Macaulayness  or the Buchsbaumness,  of $C$  in terms of $H^1_\mfr(R/I)$. 

Recall that the diameter of a $\Zset$-graded
module $M$ of finite length is  the integer $\diam(M) = \max\{n|\
M_n\neq 0\} - \min\{n|\ M_n\neq 0\} +1$ ($\diam(M) := 0$ if $M=0$).
Let $J$ be the defining ideal of an arbitrary projective curve $X$
in $\Pset^3$. Then the module  $H^1_\mfr(R/J)$ is of finite length
and let $ k(R/J)$  be the  smallest non-negative integer $k$ such
that $\mfr^k H^1_\mfr(R/J)=0$ (see \cite{MM,BCFH}).  It is obvious that  $k(R/J)\leq
\diam(H^1_\mfr(R/J))$.  
The main result of this paper states that $k(R/I) =  \diam
(H^1_\mfr(R/I))$ for an arbitrary  tetrahedral curve (see Theorem
\ref{C5}). Thus our result implies that for all
tetrahedral curves, $\diam
(H^1_\mfr(R/I))$ has no gap and $k(R/I)$ is, in this sense, as large as
possible. (Note that  monomial curves in
$\Pset^3$ also have this property, see \cite{BCFH}.) Moreover, we
can  explicitly compute $\diam (H^1_\mfr(R/I))$ in terms of
$a_1,...,a_6$ (see Theorem \ref{C2} and Theorem \ref{C5}). Since $C$ is an arithmetically Cohen-Macaulay  curve if and only if $\diam
(H^1_\mfr(R/I)) = 0$, this result  is much more general than the Francisco's one in \cite{F}. In particular, it also enables us to determine all arithmetically Buchsbaum tetrahedral curves (Theorem \ref{C7}), thus extending  Corollary 5.4 in \cite{MN}.

Our approach is to reduce the above question to a problem in
integer programming. First, based on a
description of local cohomology modules of monomial ideals given
recently in \cite{T},  we reduce the problem to describing the set of integer solutions of a certain linear constraints.   Then using the well-known Fourier-Motzkin
elimination we can determine when the set of solutions is empty (Theorem \ref{C2}). This is corresponding the case of arithmetically  Cohen-Macaulay curves. If this set is not empty, we can still use it to determine the module structure of the first  local cohomology   (Proposition \ref{C4}). Thus our result is not only   an interesting application of integer programming to Commutative Algebra, but it also shows the usefulness of Takayama's formula in \cite{T}.  We believe that  Takayama's formula, which is a generalization of Hochster's formula, can be applied in many other situations. 

The paper has four sections with the current one being an
introduction. In Section \ref{A} we recall  the main result of
Takayama in \cite{T} and relate the  problem of describing
$H^1_\mfr(R/I)$ to a problem in integer programming (Lemma
\ref{A4}). In Section \ref{B} we apply the Fourier-Motzkin
elimination to solve that integer programming problem. The structure of
the first local cohomology module is given in the last Section
\ref{C}, where the main Theorem \ref{C5} is proved and some of its
consequences are derived.

\section{Preliminaries}\label{A}

Let $I \subset R= K[x_1,...,x_n]$ be a monomial ideal. Denote by
$G(I)$ the minimal set of monomial  generators of $I$. Let $\Delta
$ be the simplicial complex corresponding to the radical ideal
$\sqrt{I}$, i.e.
$$\Delta = \{ \{i_1,...,i_k\}\subseteq \{1,...,n\}|\ x_{i_1}\cdots x_{i_k} \not\in  
\sqrt{I}\}.$$
A simplicial complex is uniquely defined by the set $\Max(\Delta
)$ of its facets.  Following \cite{T}, for $\alphaun  =
(\alpha_1,...,\alpha_n)\in \Zset^n$, we set
$$G_{\alphaun} = \{ i|\ \alpha_i <0 \},$$
and
$$\begin{array}{ll} \Delta_{\alphaun} = \{ F\subset  \{1,...,n\}\setminus  
G_{\alphaun} |\  & \text{ for\ all}\ \xund^{\betaun} = x_1^{\beta_1}\cdots  
x_n^{\beta_n} \in G(I) \  \text{there\ exists}
\\ & \ i \not\in F\cup G_{\alphaun}\ \text{such\ that}\ \beta_i > \alpha_i \geq 0\}.
\end{array}$$

\begin{Lemma}\label{A1} Denote by $I_{(x_{i_1}...x_{i_k})}$ the monomial  
ideal generated by $I$
in the localization $K[\xund]_{(x_{i_1}...x_{i_k})}$ w.r.t. the set  of all monomials in the variables
$x_{i_1},...,x_{i_k}$. Then
 $$\Delta_{\alphaun} = \{ F\subset \{1,...,n\} \setminus  G_{\alphaun} |\   
\prod_{i\not\in F \cup G_{\alphaun}}
  x_i^{\alpha_i} \not\in I_{ (\prod_{j\in F\cup G_{\alphaun}} x_j)} \}.$$
 \end{Lemma}

 \begin{pf} For simplicity we may assume that $F\cup G_{\alphaun}= \{1,...,r\}$.  
For a monomial $m\in K[x_1,...,x_n]$
 let $m' \in K[x_{r+1},...,x_n]$ be the monomial obtained from $m$ by deleting all  
powers of $x_i,\ i\leq r$.
 Let $G' = \{ m'|\ m\in G(I)\}$. Then $G'$ is a generating set of $I' :=  
I_{(x_1...x_r)}$. Note that the monomial
  $\prod_{i>r}x_i^{\alpha_i} \in I'$ if and only if there exists $m' = \prod_{i>r}  
x_i^{\beta_i} \in G'$ such
   that $\beta_i \leq \alpha_i$ for all $i>r$, or equivalently, there exists $m  
=\prod_{i=1}^n x_i^{\beta_i} \in G(I)$ such
   that $\beta_i \leq \alpha_i$ for all $i>r$. From that we immediately get the claim.
 \hfill $\square$.
 \end{pf}

 Note that all local cohomology modules $H^i_\mfr(R/I),\ i\geq 0$, inherit a natural  
$\Zset^n$-grading.  Theorem 1 in \cite{T} can be
 reformulated as follows.

 \begin{Lemma} \label{A2} Let $\rho_i = \max\{ \beta_i|\ \xund^{\betaun}\in  
G(I)\}$. For all $i\geq 0$ and $\alphaun \in \Zset^n$ we have
 $$\dim H^i_\mfr(R/I)_{\alphaun} = \begin{cases} \dim  
\tilde{H}_{i-|G_{\alphaun}| -1}(\Delta_{\alphaun},K)
 & \text{if}\ G_{\alphaun} \in \Delta \ \text{and}\ \alpha_j \leq \rho_j-1,\ j\leq n,\\
 0 &\text{otherwise}.
 \end{cases}$$
 \end{Lemma}

 From now on we consider ideals of tetrahedral curves
 $$I =(x_1,x_2)^{a_1}\cap (x_1,x_3)^{a_2}\cap (x_1,x_4)^{a_3}\cap  
(x_2,x_3)^{a_4}\cap (x_2,x_4)^{a_5}\cap (x_3,x_4)^{a_6} $$
of the polynomial ring $R=K[x_1,x_2,x_3,x_4]$.

\begin{Lemma} \label{A3} If $H^1_\mfr(R/I)_{\alphaun} \neq 0$, then $\alpha_i  
\geq 0$ for all $i\geq 1$ and
$$\Max(\Delta_{\alphaun}) =\{ \{ 1,i\}, \{j,k\}|\ \{i,j,k\} = \{2,3,4\} \}.$$
\end{Lemma}
\begin{pf} Assume $H^1_\mfr(R/I)_{\alphaun} \neq 0$. By Lemma \ref{A2},  
either $G_{\alphaun} = \emptyset $ and $\Delta_{\alphaun}$
is disconnected, or $|G_{\alphaun}| =1$ and $\Delta_{\alphaun}= \{\emptyset \}$.

If $|G_{\alphaun}| =1$, w.l.o.g. we may assume that $G_{\alphaun}=
\{1\}$, i. e. $\alpha_1 <0$ and $\alpha_2, \alpha_3, \alpha_4 \geq
0$.  By Lemma \ref{A1}, $\Delta_{\alphaun}= \{\emptyset \}$ is
equivalent to the following two conditions

(i) $x_2^{\alpha_2}x_3^{\alpha_3}x_4^{\alpha_4} \not\in I_{(x_1)} =  
(x_2,x_3)^{a_4}\cap (x_2,x_4)^{a_5}\cap (x_3,x_4)^{a_6} $, and

(ii) $x_i^{\alpha_i}x_j^{\alpha_j} \in I_{(x_1,x_k)}$ for all $\{i,j,k\} = \{2,3,4\}$.

This is impossible, because
$$ {\rm (i)} \Leftrightarrow \left[ \begin{array}{l} \alpha_2 + \alpha_3 \leq a_4-1, \  
\text{or}\\
\alpha_2 + \alpha_4 \leq a_5-1, \ \text{or}\\
\alpha_3 + \alpha_4 \leq a_6-1,
\end{array}\right.
\ \ \text{and}
\ \ {\rm (ii)} \Leftrightarrow \left\{ \begin{array}{l} \alpha_2 + \alpha_3 \geq a_4, \  
\text{and}\\
\alpha_2 + \alpha_4 \geq a_5, \ \text{and}\\
\alpha_3 + \alpha_4 \geq a_6.
\end{array}\right.$$
Hence we must have $G_{\alphaun} = \emptyset $ and
$\Delta_{\alphaun}$ is disconnected. The first  condition implies
that $\alpha_i \geq 0$ for all $i\geq 1$. Since
$\Delta_{\alphaun}$ is a disconnected simplicial complex on  a subset of  
$\{1,2,3,4\}$, in order to show the second statement of the lemma 
  it suffices to show that $\Delta_{\alphaun}$ does
not contain a facet consisting of a single point. Assume, by
contrary, that $\{1\}$ is a facet of $\Delta_{\alphaun}$. Then we
again get (i) and (ii) (the only difference now is that all
$\alpha_i \geq 0$ which, however, have no effect on (i) and (ii)).
This is  a contradiction. \hfill $\square$
\end{pf}

As an example let us consider the well-known Buchsbaum curve defined by $I= (x_1,x_2) \cap (x_3,x_4)$. In this case 
$H^1_\mfr(R/I)_{\alphaun} \neq 0$ if and only if $\alphaun = (0,0,0,0)$. We have $\Max(\Delta_{(0,0,0,0)}) =\{ \{ 1,2\}, \{3,4\}\}$.

\begin{Lemma} \label{A4} Fix an integer $d$. Assume that $\deg(\alphaun) :=  
\alpha_1 + \cdots + \alpha_4 = d$.
Then $\Max(\Delta_{\alphaun}) = \{\{1,2\}, \{3,4\}\}$ if and only if $\alphaun$  
satisfies the following system of inequalities
\begin{equation} \label{EA4}  \begin{array}{rl}
\alpha_1 + \alpha_3 & \geq a_2 \\  
\alpha_1 + \alpha_4 & \geq a_3 \\  
\alpha_2 + \alpha_3 & \geq a_4 \\ 
\alpha_2 + \alpha_4 & \geq a_5 \\
\end{array}
\end{equation}
$$ \begin{array}{rl}
\alpha_1 + \alpha_2& \leq a_1-1\\
\alpha_3 + \alpha_4 & \leq a_6-1\\
 \alpha_1 + \alpha_2 + \alpha_3 + \alpha_4 & =d\\
 \alpha_1 , \alpha_2 , \alpha_3 , \alpha_4 &\geq 0.
\end{array}
$$
In this case $\dim H^1_\mfr(R/I)_{\alphaun} =1$.
\end{Lemma}

\begin{pf} The condition $\Max(\Delta_{\alphaun}) = \{\{1,2\}, \{3,4\}\}$  implies  
$G_{\alphaun} = \emptyset$,
i.e. $\alpha_1 , \alpha_2 , \alpha_3 , \alpha_4 \geq 0$. By Lemma
\ref{A1}, $\{1,2\} \in \Delta_{\alphaun}$ if and only  if
$x_3^{\alpha_3}x_4^{\alpha_4} \not\in (x_3,x_4)^{a_6}$, or
equivalently, $\alpha_3 + \alpha_4 \leq a_6-1$. Similarly,
$\{3,4\} \in \Delta_{\alphaun}$ if and only if $\alpha_1 +
\alpha_2 \leq a_1-1$. On the other hand, $\{1,3\}, \{1,4\},
\{2,3\}, \{2,4\} \not\in \Delta_{\alphaun}$ are equivalent to the
first four inequalities given above. Thus,
$\Max(\Delta_{\alphaun}) = \{\{1,2\}, \{3,4\}\}$  implies
(\ref{EA4}). The converse is also clear from these arguments.

When $\Max(\Delta_{\alphaun}) = \{\{1,2\}, \{3,4\}\}$ we have
$\tilde{H}_0(\Delta_{\alphaun},K) \cong K$ and $|G_{\alphaun}|=0$.
Hence,  by Lemma \ref{A2}, $\dim H^1_\mfr(R/I)_{\alphaun} =1$, as
required. \hfill $\square$
\end{pf}

\section{Fourier-Motzkin elimination} \label{B}

By Lemma \ref{A4} we are interested in finding an integer solution of the  
following system of inequalities
\begin{equation} \label{EBI} \begin{array}{rl}
y_1 + y_3 & \geq a_2  \\  
y_1 + y_4 & \geq a_3 \\ 
y_2 + y_3 & \geq a_4 \\  
y_2 + y_4 & \geq a_5 \\
y_1 + y_2& \leq a_1-1\\
 y_3 + y_4 & \leq a_6-1\\
 y_1 + y_2 + y_3 + y_4 & =d\\
 y_1 , y_2 , y_3 , y_4 &\geq 0.
\end{array}
\end{equation}
For this purpose we apply the Fourier-Motzkin elimination which at
first enables to find a real solution  of a system of linear
equalities and inequalities, see, e.g. \cite{DT}, Section 2.3. We
sketch here the algorithm by considering a concrete example.

{\it Example}.  Consider the system
\begin{equation} \label{EB1}  \begin{array}{rl}
y_1 + 2y_2 - y_3 + 4 & \geq 0\\
- 2y_1 + y_2 + 3y_3 - 2&\geq 0\\
2y_2 - y_3 &\geq 0\\
y_1 &= y_2+y_3.
\end{array}
\end{equation}
First, replace the equality $y_1 = y_2+y_3$ by two inequalities
$y_1 \geq  y_2+y_3$ and $y_1 \leq  y_2+y_3$.  The obtained system
is not reduced w.r.t. $y_1$, i.e. $y_1$ appears with a non-zero
coefficient in at least one inequality. After dividing by the
absolute value of the coefficient of $y_1$  when nonzero and
rearranging the terms and the order of the constraints, we can
then partition them  in 3 groups, depending on whether in a
particular constraint $y_1$ is on the right or the left hand, or its $y_1$-coefficient  
is zero.
$$ \begin{array}{rlr}
\frac{1}{2} y_2 + \frac{3}{2} y_3 - 1 & \geq y_1 & \hskip1cm (E1)\\
y_2 + y_3 &\geq y_1 & (E2)\\
y_1 &\geq  - 2y_2+ y_3 - 4  & \ (E3)\\
y_1 &\geq  y_2 + y_3  & (E4)\\
2y_2 - y_3 &\geq 0. &(E5)
\end{array}$$
Combining each inequality  in the first group  $\{ (E1),\ (E2)\}$
with another one  in the second group  $\{ (E3),\ (E4) \}$  and keep
all inequalities  in the third group ($\{(E5) \}$ in this example), we
obtain a new system of inequalities
$$ \begin{array}{rlr}
\frac{1}{2} y_2 + \frac{3}{2} y_3 - 1 & \geq  - 2y_2+ y_3 - 4  & \hskip1cm  
(E1,E3)\\
\frac{1}{2} y_2 + \frac{3}{2} y_3 - 1 & \geq  y_2 + y_3  & (E1,E4)\\
y_2 + y_3 &\geq - 2y_2+ y_3 - 4  &  (E2,E3)\\
y_2 + y_3 &\geq  y_2 + y_3  & (E2,E4)\\
2y_2 - y_3 &\geq 0. &(E5)
\end{array}$$
The temporary label $(E1,E3)$ means that this inequality appears
by combining $(E1)$ and $(E3)$. Note that in  the last system,
$(E1,E3)$ follows from $(E1,E4)$ and $(E2,E3)$.  For short, we
will write this reduction as $(E1,E4) +(E2,E3)\Rightarrow
(E1,E3)$. The  constraint  $(E2,E4)$ trivially holds. We say that
$(E1,E3)$ and $(E2,E4)$ are redundant.  Deleting the redundant
inequalities, we finally get the system
\begin{equation} \label{EB2}  \begin{array}{rl}
\frac{1}{2} y_2 + \frac{3}{2} y_3 - 1 & \geq  y_2 + y_3  \\
y_2 + y_3 &\geq - 2y_2+ y_3 - 4  \\
2y_2 - y_3 &\geq 0.
\end{array}
\end{equation}
Thus  (\ref{EB1}) implies (\ref{EB2}), where $y_1$ appears with
zero coefficient in all inequalities.  We say that $y_1$ has been
"eliminated". The process is repeated with the new system except
now $y_2$ is eliminated.

We now apply the Fourier-Motzkin elimination to our system (\ref{EBI}). First  
rewrite it in the form
\begin{equation} \label{EBII} \begin{array}{rl}
a_6-1 -y_3 &\geq y_4 \\  
y_4 &= d -y_1 - y_2 - y_3 \\ 
y_4 & \geq a_3 -y_1 \\  
y_4  &\geq a_5 -y_2  \\  
y_4 &\geq 0\\
y_1 + y_3 & \geq a_2\\
 y_1 + y_2& \leq a_1-1\\
 y_2 + y_3 & \geq a_4\\
 y_1 , y_2 , y_3 &\geq 0.
\end{array}
\end{equation}
Eliminating $y_4$ we then get
\begin{equation} \label{EBIII}  \begin{array}{rl}
d-a_3-y_2 &\geq y_3 \\ 
d-a_5-y_1& \geq y_3  \\  
d-y_1-y_2 & \geq y_3 \\  
y_3&\geq 0\\
y_3 & \geq a_2-y_1\\
y_3 & \geq a_4-y_2\\
 y_1+y_2+a_6-d-1 &\geq 0\\
 y_1+y_2 &\leq a_1-1\\
 y_1,y_2&\geq 0.
\end{array}
\end{equation}
Eliminating $y_3$  we now obtain
\newpage
\begin{equation} \label{EBIV}    \begin{array}{rlr}
d-a_3 &\geq y_2 &\hskip1cm (\ref{EBIV}.1) \\
d-a_2-a_3+y_1 &\geq y_2 & (\ref{EBIV}.2)\\
d-y_1 &\geq y_2 & (\ref{EBIV}.3)\\
d-a_2 &\geq y_2 & (\ref{EBIV}.4)\\
a_1-1-y_1 &\geq y_2 & (\ref{EBIV}.5)\\
y_2 &\geq 0 & (\ref{EBIV}.6)\\
y_2 &\geq a_4+a_5 - d+y_1 & (\ref{EBIV}.7)\\
y_2 &\geq d+1-a_6 -y_1 & (\ref{EBIV}.8)\\
d-a_5 &\geq y_1 &\hskip1cm (\ref{EBIV}.9)\\
d-a_4 & \geq y_1 & (\ref{EBIV}.10)\\
y_1 &\geq 0 & (\ref{EBIV}.11)\\
d & \geq a_3+a_4 & (\ref{EBIV}.12)\\
d & \geq a_2+a_5. & (\ref{EBIV}.13)
\end{array}
\end{equation}
By eliminating $y_2$ we get a system of 20 constraints. However 7
of them  are redundant:  $(\ref{EBIV}.12) \Rightarrow
(\ref{EBIV}.1, \ref{EBIV}.6);\ (\ref{EBIV}.9) + (\ref{EBIV}.12)
\Rightarrow (\ref{EBIV}.1, \ref{EBIV}.7);\
(\ref{EBIV}.12)+(\ref{EBIV}.13) \Rightarrow
(\ref{EBIV}.2,\ref{EBIV}.7); \ (\ref{EBIV}.9)+(\ref{EBIV}.10)
\Rightarrow (\ref{EBIV}.3,\ref{EBIV}.6),
(\ref{EBIV}.3,\ref{EBIV}.7), (\ref{EBIV}.4,\ref{EBIV}.7)$ and
$(\ref{EBIV}.13) \Rightarrow (\ref{EBIV}.4,\ref{EBIV}.6)$. Deleting
these redundant  constraints we get
\begin{equation} \label{EBV}  \begin{array}{rlr}
d-a_4 &\geq y_1 &\hskip1cm (\ref{EBV}.1) \\
d-a_5 &\geq y_1 & (\ref{EBV}.2)\\
a_1-1 &\geq y_1 & (\ref{EBV}.3)\\
   \lfloor  \frac{1}{2}(d+a_1-a_4-a_5-1) \rfloor&\geq y_1 & (\ref{EBV}.4)\\
y_1 &\geq 0 & (\ref{EBV}.5)\\
y_1 &\geq \lceil \frac{1}{2}(a_2+a_3-a_6+1)\rceil & (\ref{EBV}.6)\\
y_1 &\geq a_2+a_3 -d & (\ref{EBV}.7)\\
y_1 & \geq a_3 -a_6 +1 & (\ref{EBV}.8)\\
y_1 &\geq a_2-a_6 + 1 & (\ref{EBV}.9)\\
a_1+a_6-2 &\geq d  & (\ref{EBV}.10)\\
d &\geq a_2+a_5 & (\ref{EBV}.11)\\
\end{array}
\end{equation}
$$\begin{array}{rlr}
d &\geq a_3 + a_4 & \hskip1cm  (\ref{EBV}.12)\\
a_6&\geq 1. & (\ref{EBV}.13)
\end{array}$$
Here, for a real number $a$, we set 
$$\lceil a \rceil = \min\{ n\in
\Zset|\ n\geq a\}\ \text{and}\ \lfloor a \rfloor = \max\{ n\in \Zset|\
n\leq a\}.$$
 Eliminating $y_1$ we get a system of 24 constraints.
Among them 14 are redundant: $(\ref{EBV}.12) \Rightarrow (\ref{EBV}.1,
\ref{EBV}.5);\ (\ref{EBV}.1,\ref{EBV}.9)+ (\ref{EBV}.12)
\Rightarrow (\ref{EBV}.1,\ref{EBV}.6);\ (\ref{EBV}.11) +
(\ref{EBV}.12) \Rightarrow  (\ref{EBV}.1, \ref{EBV}.7);\
(\ref{EBV}.12) + \ref{EBV}.13) \Rightarrow
(\ref{EBV}.1,\ref{EBV}.8);\ (\ref{EBV}.11) \Rightarrow
(\ref{EBV}.2,\ref{EBV}.5);\ (\ref{EBV}.2,\ref{EBV}.8) +
(\ref{EBV}.11) \Rightarrow (\ref{EBV}.2,\ref{EBV}.6);\
(\ref{EBV}.11) + (\ref{EBV}.12) \Rightarrow  (\ref{EBV}.2,
\ref{EBV}.7);\ (\ref{EBV}.11) + (\ref{EBV}.13) \Rightarrow
(\ref{EBV}.2, \ref{EBV}.9);\
(\ref{EBV}.3,\ref{EBV}.7)+(\ref{EBV}.10) \Rightarrow
(\ref{EBV}.3,\ref{EBV}.6);\ (\ref{EBV}.10) +
(\ref{EBV}.12)\Rightarrow (\ref{EBV}.3,\ref{EBV}.8);$ 
$ (\ref{EBV}.10)+(\ref{EBV}.11) \Rightarrow
(\ref{EBV}.3,\ref{EBV}.9);\ 
(\ref{EBV}.11)+(\ref{EBV}.12)+(\ref{EBV}.3,\ref{EBV}.7)
\Rightarrow (\ref{EBV}.4,\ref{EBV}.7);\
(\ref{EBV}.10)+(\ref{EBV}.12)+(\ref{EBV}.2,\ref{EBV}.8)
\Rightarrow (\ref{EBV}.4,\ref{EBV}.8)$  and $
(\ref{EBV}.10)+(\ref{EBV}.11)+(\ref{EBV}.1,\ref{EBV}.9)
\Rightarrow (\ref{EBV}.4,\ref{EBV}.9).$ Deleting these redundant  constraints, we  
finally
get the system
\begin{equation} \label{EBVI}  \begin{array}{rl}
a_1+a_6-2 &\geq d \\
d &\geq a_2+a_5 \\  
d&\geq a_3+a_4\\   
d &\geq a_2+a_4-a_6 +1\\
d& \geq a_3+a_5-a_6 +1\\
d& \geq a_2+a_3-a_1+1\\
d&\geq a_4+a_5-a_1+1 \\
 \lfloor  \frac{1}{2}(d+a_1-a_4-a_5-1)\rfloor &\geq \lceil  
\frac{1}{2}(a_2+a_3-a_6+1)\rceil \\
a_1,a_6 &\geq 1.
\end{array}
\end{equation}

\begin{Lemma} \label{B1} Assume that $a_1+a_6 - 2\geq d \geq \max\{a_2+a_5,  
a_3+a_4\}$. Then $\lfloor
\frac{1}{2}(d+a_1-a_4-a_5-1)\rfloor < \lceil
\frac{1}{2}(a_2+a_3-a_6+1)\rceil$ if and only if $a_2+a_3-a_6$ is
even and $a_1+a_6-2 = a_2+a_5 = a_3+a_4$.
\end{Lemma}

\begin{pf} If $a_2+a_3-a_6$ is odd, then
$$\lceil \frac{1}{2}(a_2+a_3-a_6+1)\rceil = \frac{1}{2}(a_2+a_3-a_6+1).$$
Since $a_2+a_5+a_3+a_4 \leq d + a_1+a_6 -2$, we get $d+a_1-a_4-a_5
-1 \geq a_2+a_3-a_6+1$, which yields
$$\frac{1}{2}(d+a_1-a_4-a_5-1) \geq \frac{1}{2}(a_2+a_3-a_6+1).$$
Hence
$$\lfloor \frac{1}{2}(d+a_1-a_4-a_5-1)\rfloor \geq \frac{1}{2}(a_2+a_3-a_6+1) =  
\lceil \frac{1}{2}(a_2+a_3-a_6+1)\rceil.$$
If $a_2+a_3-a_6$ is even, then
$$\lceil \frac{1}{2}(a_2+a_3-a_6+1)\rceil = \frac{1}{2}(a_2+a_3-a_6) + 1.$$
In the case  $a_1+a_6-2 > \min\{ a_2+a_5, a_3+a_4\}$, we have  
$a_2+a_5+a_3+a_4
\leq d + a_1+a_6 -3$. Hence $d+a_1-a_4-a_5 -1 \geq a_2+a_3-a_6+2$,
which implies
$$\lfloor \frac{1}{2}(d+a_1-a_4-a_5-1)\rfloor \geq \frac{1}{2}(a_2+a_3-a_6)  + 1=  
\lceil \frac{1}{2}(a_2+a_3-a_6+1)\rceil.$$
The left case is $a_1+a_6-2 = \min\{ a_2+a_5, a_3+a_4\}$. Since $a_1+a_6-2 \geq  
d\geq  \max\{ a_2+a_5, a_3+a_4\}$, we must have $d= a_2+a_5 = a_3+a_4 =  
a_1+a_6-2$. Then
$d+a_1-a_4-a_5-1 = a_2+a_3-a_6+1$ is an odd number. Therefore
$$\lfloor \frac{1}{2}(d+a_1-a_4-a_5-1)\rfloor <  \lceil  
\frac{1}{2}(a_2+a_3-a_6+1)\rceil.$$
This completes the proof of  the lemma. \hfill $\square$
\end{pf}

Going back from (\ref{EBVI}) to (\ref{EBII}), the Fourier-Motzkin
algorithm gives us in general only a rational solution  of
(\ref{EBI}) if (\ref{EBVI}) holds. However, in our concrete
situation we can already find an integer solution.

\begin{Lemma} \label{B2} Let
$$\begin{array}{ll}
\Acal = \max \{ & a_2 + a_5, \ a_3+a_4, \ a_2+a_4-a_6+1, \ a_3+a_5-a_6+1, \\
& a_2+a_3-a_1+1, \ a_4+a_5-a_1+1\}.
\end{array}$$
The system (\ref{EBI}) has an integer solution if and only if $a_1,a_6 \geq 1$ and  
one of the following conditions holds:
\begin{itemize}
\item[\rm{(i)}] $a_1+a_6 - 2>\Acal$ and $a_1+a_6 - 2\geq d \geq \Acal$.
\item[\rm{(ii)}] $a_1+a_6 - 2=\Acal =d$ and $a_1+a_6 - 2 > \min\{a_2+a_5,\  
a_3+a_4\}$.
\item[\rm{(iii)}] $a_1+a_6 - 2=a_2+a_5 = a_3+a_4 = \Acal= d$ and $a_2+a_3-a_6$ is  
odd.
\end{itemize}
\end{Lemma}

\begin{pf} If (\ref{EBI}) has an integer   solution, then by Fourier-Motzkin  
algorithm, (\ref{EBVI}) holds. Using Lemma \ref{B1} we get the necessity.

Assume that $a_1,a_6 \geq 1$ and one of the above conditions
(i)-(iii) holds. Then for any $d$ such that $\Acal \leq d \leq a_1+a_6-2$, the
system (\ref{EBVI}) holds by Lemma \ref{B1}. Fix such an integer
$d$. Denote by $\Lcal_{\ref{EBV}}$ the minimum of integers in the left
sides of $(\ref{EBV}.1)-(\ref{EBV}.4)$ and $\Rcal_{\ref{EBV}}$ the
maximum of integers in the right sides of
$(\ref{EBV}.5)-(\ref{EBV}.9)$. Then from (\ref{EBVI}) it follows
that $\Lcal_{\ref{EBV}} \geq \Rcal_{\ref{EBV}}$. Hence $y_1 =
\Rcal_{\ref{EBV}}$ is an integer solution of (\ref{EBV}). Putting $y_1
= \Rcal_{\ref{EBV}}$ into (\ref{EBIV})-(\ref{EBII}) and repeating this
process, we can similarly define $\Lcal_{\ref{EBIV}} \geq
\Rcal_{\ref{EBIV}},\  \Lcal_{\ref{EBIII}} \geq \Rcal_{\ref{EBIII}},\
\Lcal_{\ref{EBII}} \geq \Rcal_{\ref{EBII}}$ such that $y_1 =
\Rcal_{\ref{EBV}},\ y_2 = \Rcal_{\ref{EBIV}},\ y_3 = \Rcal_{\ref{EBIII}},\ y_4
 = \Rcal_{\ref{EBII}}$ is an integer solution of (\ref{EBII}),
which is equivalent to (\ref{EBI}) . \hfill $\square$
\end{pf}

\section{Structure of the first local cohomology module} \label{C}

In this section we describe the first local cohomology module of $R/I$. From now  
on, w.l.o.g.,  we always assume that $a_1+a_6$ is the maximum among the sums $a_1+a_6,\ a_2+a_5,\ a_3+a_4$. In other words we may assume that the following holds:
$$ (*) \hskip2cm a_1 + a_6 \geq \max\{ a_2+a_5, \ a_3+a_4\}.$$

\begin{Lemma} \label{C1} Under the assumption (*) there exists no $\alphaun \in  
\Zset^4$ such that
$\Max(\Delta_{\alphaun}) = \{ \{1,3\},\ \{2,4\}\}$ or $\Max(\Delta_{\alphaun}) =  
\{ \{1,4\},\ \{2,3\}\}$.
\end{Lemma}

\begin{pf} Assume, w.l.o.g., the existence of $\alphaun \in \Zset^n$ such that  
$\Max(\Delta_{\alphaun}) = \{ \{1,3\},\ \{2,4\}\}$. Then
applying Lemma \ref{A4} and Lemma \ref{B2} to this situation we would get  
$a_2+a_5-2 \geq a_1+a_6$, a contradiction to (*). \hfill $\square$
\end{pf}

We can now  explicitly determine all arithmetically Cohen-Macaulay
tetrahedral  curves in terms of $a_i$. This result recovers the
main theorem in \cite{F}.

\begin{Theorem} \label{C2} Let 
$$\begin{array}{ll}
\Acal = \max \{ & a_2 + a_5, \ a_3+a_4, \ a_2+a_4-a_6+1, \ a_3+a_5-a_6+1, \\
& a_2+a_3-a_1+1, \ a_4+a_5-a_1+1\}.
\end{array}$$
Under the assumption (*), a  tetrahedral
curve $C(a_1,...,a_6)$ is arithmetically Cohen-Macaulay if and only if  one of the  
following conditions holds:
\begin{itemize}
\item[\rm{(i)}] $a_1=0$ or $a_6 =0$;
\item[\rm{(ii)}] $a_1+a_6 - 2<\Acal $;
\item[\rm{(iii)}] $a_1+a_6 - 2=a_2+a_5 = a_3+a_4 = \Acal$ and $a_2+a_3-a_6$ is  
even.
\end{itemize}
\end{Theorem}

\begin{pf} By Lemma \ref{C1}, $C= C(a_1,...,a_6)$ is arithmetically  
Cohen-Macaulay if and only
if there is no $d$ such that the system (\ref{EBI}) has an integer solution. Hence  
the statement follows from Lemma \ref{B2}. \hfill $\square$
\end{pf}

 {\bf Remark}. In \cite{MN}, Question 7.4(5), Migliore and Nagel asked whether  
an arithmetically  Cohen-Macaulay
 tetrahedral curve $C= C(a_1,...,a_6)$ can be explicitly identified by the 6-tuples  
$a_1,...,a_6$. This question was solved
 by Francisco in \cite{F}. His main result says that under the assumption (*),  
$C(a_1,...,a_6)$ is arithmetically Cohen-Macaulay if and
 only if  one of the following conditions holds:
 \begin{itemize}
\item[\rm{(a)}] $a_1=0$ or $a_6 =0$;
\item[\rm{(b)}] $a_1+a_6 = \epsilon + \max\{a_2+a_5,\ a_3+a_4\}$, where  
$\epsilon \in \{0,1\}$.
\item[\rm{(c)}] $2a_1 < a_2 + a_3 -a_6 +3$ or $2a_1 < a_4 + a_5 -a_6 +3$ or  
$2a_6< a_2 + a_4 -a_1 +3$ or $2a_6 < a_3 + a_5 -a_1 +3$;
\item[\rm{(d)}] All inequalities of (c) fail, $a_1+a_6 =a_2+a_5 +2= a_3+a_4 +2$  
and $a_1+a_3+a_5$ is even.
\end{itemize}
 One can easily check that this statement is equivalent to that of Theorem \ref{C2}.

 Assume now that $C$ is not arithmetically Cohen-Macaulay. Then $a_1,a_6\geq  
1$ and one of three conditions in Lemma \ref{B2}
 is satisfied. In particular $\Acal \leq a_1+a_6-2$. Let
$$\begin{array}{ll} T_1  & =\{ \yund \in \Nset^4|\   y_1 + y_3  \geq a_2  ,\ y_1 +  
y_4  \geq a_3 ,\  y_2 + y_3  \geq a_4 ,\  y_2 + y_4  \geq a_5 \},\\
T_2 & = \{ \yund \in T_1|\ y_1 + y_2 \geq a_1\},\end{array}$$
and
$$T_3 = \{ \yund \in T_1|\ y_3 + y_4 \geq a_6\}.$$
Let $S= T_1\setminus (T_2\cup T_3)$. Then the set $S_d$ of all
elements of degree $d$  of $S$ is the set of all solutions  of the system
(\ref{EBI}).  As usual we identify $K[T_i],\ i\leq 3$,  and $K[S]$
with  subsets of $R=K[x_1,...,x_4]$. Note that $K[T_i], \ i\leq
3,$ are ideals  of  $R$. Hence we may consider $K[S]$ as a factor
module $K[T_1]/K[T_2]+K[T_3]$. Thus, the  module structure on
$K[S]$ over $R$ is defined as follows: for $\alphaun \in S$ and
$\betaun \in \Nset^4$,
 $$\xund^{\betaun}\cdot \xund^{\alphaun} = \begin{cases} \xund^{\betaun+  
\alphaun} \ &\text{if} \ \betaun+ \alphaun \in S,\\
 0 \ & \text{otherwise}.\end{cases}$$

The following result describes the module structure of $H^1_\mfr(R/I)$.

 \begin{Proposition}\label{C4} Under the assumption (*), $$H^1_\mfr(R/I) \cong  
K[S]$$ as graded modules over $R$.
 \end{Proposition}

 \begin{pf} Let
$${\mathcal C}^\bullet : 0 \rightarrow R/I \rightarrow \oplus_{i=1}^4 (R/I)_{x_i}   
\rightarrow \cdots \rightarrow (R/I)_{x_1x_2x_3 x_4}
 \rightarrow 0,$$
be the $\check{C}$ech complex of $R/I$. Then $H^1_\mfr(R/I) \cong
H^1({\mathcal C}^\bullet)$. By \cite{T}, Lemma 2, for all
$\alphaun\in \Zset^4$ there is an isomorphism of complexes
$$({\mathcal C}^\bullet_{\alphaun}) \cong \Hom_{\Zset}({\mathcal  
C}(\Delta_{\alphaun})[-j-1],K),$$
where $j= |G_{\alphaun}|$ and ${\mathcal
C}(\Delta_{\alphaun})[-j-1]$ means the  shifting  of the augmented
oriented chain complex ${\mathcal C}(\Delta_{\alphaun})$ by $-j-1$.
Denote by $\pi $ the simplicial complex on $\{1,2,3,4\}$ with
$\Max(\pi) = \{ \{1,2\},\ \{3,4\}\}$. By Lemmas
\ref{A3}, \ref{A4} and \ref{C1} it follows that $H^1({\mathcal
C}^\bullet_{\alphaun}) \neq 0$ if and only if $\Delta_{\alphaun} =
\pi $, $G_{\alphaun} =\emptyset $ and $\alphaun \in S$. Moreover,
in this case $H^1({\mathcal C}^\bullet_{\alphaun}) \cong
K\xund^{\alphaun}$. From this we get $H^1_\mfr(R/I) \cong K[S]$,
as required. 

\rightline{ $\square$}
\end{pf}

The above description of $S$ allows us  to describe the module structure of $K[S]$ in an obvious way. Of course, $S$ can be 
written as:
$$\begin{array}{ll} S= \{ \yund \in \Nset^4|\  & y_1 + y_3  \geq a_2  ,\ y_1 +  
y_4  \geq a_3 ,\ y_2 + y_3  \geq a_4 ,\  y_2 + y_4  \geq a_5, \\
& y_1+y_2 < a_1,\ y_3+y_4 < a_6 \}.
\end{array}$$
It is easy to write a program to compute this set $S$. Hence the module structure of $H^1_\mfr(R/I) $ is known once $a_1,...,a_6$ are given.

We say that a non-zero $\Zset$-graded module $M$ has no gap if
$M_i \neq 0$ and $M_j\neq 0$ for  some $i\leq j$, then $M_k\neq 0$
for all $i\leq k\leq j$. Recall that the diameter of a module $M$
of finite length is defined as
$$\diam (M) = \ed(M) - \bg(M) +1,$$
where $\bg(M) =\min\{i|\ M_i \neq 0\}$ and $\ed(M) =\max\{i|\ M_i \neq 0\}$  (if  
$M=0$ we set $\diam(M) =0$).

\begin{Theorem} \label{C5} Let
$$\begin{array}{ll}
\Acal = \max \{ & a_2 + a_5, \ a_3+a_4, \ a_2+a_4-a_6+1, \ a_3+a_5-a_6+1, \\
& a_2+a_3-a_1+1, \ a_4+a_5-a_1+1\}.
\end{array}$$
Assume that (*) holds and the tetrahedral curve $C$  
is not arithmetically Cohen-Macaulay. Then
$a_1+a_6-2 \geq \Acal$ and
$$k(R/I) =  \diam (H^1_\mfr(R/I)) = a_1+a_6- \Acal -1.$$
In particular, $H^1_\mfr(R/I)$ has no  
gap.
\end{Theorem}

\begin{pf} Since $R/I$ is not a Cohen-Macaulay ring, by Theorem \ref{C2},  
$a_1+a_6-2 \geq \Acal$ and $a_1,a_6\geq 1$.
By Lemma \ref{B2}, for each $d$ such that $\Acal \leq d \leq a_1+a_6-2$ we have $S_d \neq  
\emptyset $. Hence,
 by Proposition \ref{C4},  $H^1_\mfr(R/I)$ has no gap, $\bg(H^1_\mfr(R/I)) = \Acal$ and  
$\ed(H^1_\mfr(R/I)) = a_1+a_6-2$,
 which implies $\diam (H^1_\mfr(R/I)) = a_1+a_6- \Acal -1$.

Further, let $\alphaun = (\alpha_1,\alpha_2,\alpha_3,\alpha_4)\in
S_\Acal$ be a fixed element. Then  $\alpha_1+ \alpha_2 \leq a_1-1$ and
$\alpha_3+\alpha_4 \leq a_6-1$. Let $\alphaun^* = (\alpha_1, a_1-1
-\alpha_1,\alpha_3, a_6-1-\alpha_3)$. Since $a_1-1 -\alpha_1\geq
\alpha_2$ and $ a_6-1-\alpha_3\geq \alpha_4$, the condition
$\alphaun\in T_1$ implies $\alphaun^* \in T_1$ too. On the other
hand $\alphaun^* \not\in T_1\cup T_2$. Hence $\alphaun^* \in
S_{a_1+a_6-2}$. Note that $\alphaun^* = \alphaun + \betaun$, where
$\betaun= (0, a_1-1 -\alpha_1-\alpha_2, 0, a_6-1-\alpha_3-
\alpha_4)\in \Nset^4$ and $\deg (\betaun) = a_1+a_6-\Acal -2$.
Therefore, by Proposition \ref{C4}, $$\xund^{\betaun}
H^1_\mfr(R/I)_{\alphaun} \cong H^1_\mfr(R/I)_{\alphaun + \betaun}
= H^1_\mfr(R/I)_{\alphaun^*} \neq 0,$$ which yields
$$k(R/I) \geq a_1+a_6- \Acal-1 = \diam (H^1_\mfr(R/I)).$$
Since  $\diam (H^1_\mfr(R/I)) \geq k(R/I)$, we finally get $k(R/I)  = \diam  
(H^1_\mfr(R/I))$, as required. \hfill $\square$
\end{pf}

In the above proof we already showed:

\begin{Corollary} Assume that (*) holds and the tetrahedral curve $C$ is not  
arithmetically Cohen-Macaulay. Then
$a_1+a_6-2 \geq \Acal$ and $\ed(H^1_\mfr(R/I)) = a_1+a_6-2$.
\end{Corollary}

Recall that $C$ is  arithmetically Buchsbaum if and only if
$k(R/I) \leq1$. As an immediate consequence of  Theorem \ref{C5}
we recover Corollary 4 in \cite{MN}.

\begin{Corollary} A tetrahedral curve $C$ is  arithmetically Buchsbaum if and  
only if
$$H^1_\mfr(R/I) \cong K^m(t),$$
for some non-negative integers $m,t$.
\end{Corollary}

 Migliore and Nagel found all arithmetically Buchsbaum tetrahedral curves which  
are so-called minimal (see Corollary \ref{C8} below).
  Using Theorem \ref{C5} and  \ref{C2} we are able to determine all arithmetically  
Buchsbaum tetrahedral curves which are not necessarily minimal.

 \begin{Theorem} \label{C7} Under the assumption (*), a tetrahedral curve $C$ is   
 arithmetically Buchsbaum if and only if   one of the following conditions is satisfied:  

(i) $a_1=0$ or $a_2=0$;

(ii) $a_1+a_6 - 2 \leq \Acal $.
\end{Theorem}

\begin{pf} If $C$ is arithmetically Cohen-Macaulay, by Theorem \ref{C2},  one of the above condition holds. Assume that $C$ is  not arithmetically Cohen-Macaulay and arithmetically  
Buchsbaum. Then $k(R/I) =1$.
By Theorem \ref{C5}, $a_1,a_6\geq 1$ and $a_1+a_6 - 2=\Acal $. Conversely, by Theorem \ref{C2} we may assume from the beginning that $a_1,a_6\geq 1$.  Under  
these conditions,  again by Theorem \ref{C5},we immediately
 have $k(R/I) \leq 1$, i.e. $C$ is  arithmetically Buchsbaum.  \hfill $\square$
\end{pf}

Migliore and Nagel introduced the following notion: Assume that
\newline
$a_6= \max\{a_1,...,a_6\}$. A tetrahedral  curve $C$ is said to be
{\it minimal} if $a_1> \max\{ a_2+a_4,\ a_3+a_5\}$ and $a_6>
\max\{ a_2+a_3,\ a_4+a_5\}$ (see \cite{MN}, Definition 3.4 and
Corollary 3.5). Note that in this case we already have $a_1,a_6
\geq 1$ and $a_1+a_6 - 2\geq \Acal$.

\begin{Corollary} {\rm (\cite{MN}, Corollary 4.3 and Corollary 5.4)}. \label{C8}  
Assume that 
\newline $a_6= \max\{a_1,...,a_6\}$ and $C$ is a
minimal tetrahedral curve. Then
\begin{itemize}
\item[{\rm  (i)}] $C$ is not arithmetically Cohen-Macaulay.
\item[{\rm  (ii)}] $C$ is arithmetically Buchsbaum if and only if
either $a_2=a_5=0$ and $a_1=a_6 = a_3+1=a_4+1$ or $a_3=a_4=0$  and
$a_1=a_6 = a_2+1=a_5+1$.
\end{itemize}
\end{Corollary}

\begin{pf} Since $a_1> \max\{ a_2+a_4,\ a_3+a_5\}$ and $a_6> \max\{  
a_2+a_3,\ a_4+a_5\}$, we have
\begin{equation} \label{EC8} \begin{array}{rl}
a_1+a_6-2 \geq \max\{ &  a_2 + a_5 + 2a_4,\ a_2+a_5 +2a_3,\\
& a_3+a_4+2a_2, \ a_3+a_4+2a_5\} \geq \Acal.
\end{array}
\end{equation}
If $C$ is arithmetically Buchsbaum, then since $a_1,a_6\geq 1$, by
Theorem \ref{C2} and Theorem \ref{C7}, we must have $a_1+a_6-2
=\Acal$.  Combining with (\ref{EC8}) it implies that either
$a_2=a_5=0$ or $a_3=a_4=0$. W.l.o.g. assume that $a_2=a_5=0$. Then
$\Acal = a_3+a_4$ and $a_1+a_6-2= a_3+a_4$. Since $a_1, a_6>\max\{
a_3,\ a_4\}$, the later equality gives $a_1=a_6 = a_3+1=a_4+1$. In
this case $a_2+a_3-a_6 = -1$ is odd, so $C$ is not arithmetically
Cohen-Macaulay. Thus we have proved (i) and the necessity of (ii).
The sufficiency of (ii) immediately follows from Theorem \ref{C7}.
\hfill $\square$
\end{pf}

Similarly, using Theorem \ref{C5}, we can quickly get

\begin{Corollary} {\rm (\cite{MN}, Lemma 6.2)}. \label{C9} Assume that $a_6=  
\max\{a_1,...,a_6\}$ and $C$ is
a minimal tetrahedral curve. Then $\diam (H^1_\mfr(R/I)) = 2$ if and only if after a  
suitable permutation of variables we have
$(a_1,...,a_6) = (k,k-1,0,0,k-1,k+1),\ k\geq 1$ or $(a_1,...,a_6) = (k,k-2,0,0,k-1,k),\  
k\geq 2$.
\end{Corollary}

\end{document}